\documentclass{article}
 \usepackage[ansinew]{inputenc}
 \usepackage[psamsfonts]{amssymb}
  \def\bk{{l\!k}}

\title{Torsion primes in loop space homology}
\author{Y. Felix, S. Halperin and J.-C.  Thomas}
\begin{document}
\maketitle

\begin{abstract}
We construct a finite 1-connected CW complex $X$ such that
$H_*(\Omega X;\mathbb Z)$ has $p$-torsion for the infinitely many
primes satisfying $p\equiv 5,7,17,19$ mod $24$, but no $p$-torsion
for the infinitely many primes satisfying $p\equiv 13$ or $23$ mod
$24$.
\end{abstract}

 \vspace{5mm}\noindent {\bf AMS Classification} : 55P35

 \vspace{2mm}\noindent {\bf Key words} : loop space homology

\vspace{5mm}

A   {\sl homology torsion prime} for a topological space $Y$ is a prime $p$ such that
  $H_*(Y;\mathbb Z)$ has $p$ torsion. In \cite{An} and \cite{Av} respectively, D. Anick
 and L. Avramov   constructed simply
connected finite CW complexes $X$ of dimension 4 such that $H_*(\Omega
X,\mathbb Z)$ has   $p$-torsion   for all primes $p$.
Since there are only countably many   homotopy types of simply connected finite CW
complexes, this immediately raised the

\vspace{2mm}\noindent {\bf Question.} If $\cal P$ denotes the set of all primes, what are the conditions for
a subset ${\cal D} \subset \cal P$, to be the set of homology torsion primes of
the loop space   of a finite simply connected CW complex ?

\vspace{2mm} Of course, any finite set can be realized. It suffices to take a product of Moore spaces.
On the other hand by a result of McGibbon and Wilkerson \cite{MW},
if $X$ is rationally elliptic (i.e., dim $\pi_*(X)\otimes \mathbb Q <\infty$),
then the set of homology torsion primes for $\Omega X$
is  finite.

 In this note we adapt the construction of Anick to prove

\vspace{3mm}\noindent {\bf Theorem 1.} {\sl  There exists
  a 4-dimensional finite simply connected  CW complex $X$
  such that the homology torsion primes for
  $\Omega X$ are an infinite set  in $\cal P$, with an infinite complement.}

\vspace{3mm}\noindent {\bf Theorem 2.} {\sl  For each finite subset $E$ of   ${\cal P} $
 there exists
  a 4-dimensional finite simply connected  CW complex $X_E$  which has ${\cal P}\backslash E$
  as the set of homology torsion primes for $\Omega X_E$.}

\vspace{3mm}

Our examples are CW complexes $X$ obtained by attaching thirteen 4-cells to
a wedge of eight $2$-spheres
via a map $f : Y = \vee_{i=1}^{13} S^3_i \to Z=\vee_{j=1}^8 S_j^2$.
   We label the 2-spheres   $u_1, u_2, u_3, u_4, v,w, x_1,
x_2$ and, fixing six integers $a,b,c,d,a_2, b_2$, attach the 13 cells   in degree 4   along the elements
$$
\begin{array}{c}
[x_1 ,u_1] - [u_1,v] - a[u_2,v] \,, \hspace{5mm}
 [x_1,u_2] - b[u_2,v] - d [u_3,v]\,, \hspace{5mm}
 [x_1 , u_3] - c[u_2,v]\\
 \mbox{} [x_1,u_4]\,, \hspace{5mm}  [x_1,v] \,, \hspace{5mm} [x_1,w]
\,, \hspace{5mm}
 [x_2,u_1]\,, \hspace{5mm}  [x_2,u_3] \,, \hspace{5mm} [x_2,u_4]\\
 \mbox{}[x_2,v]\,, \hspace{5mm} [x_2,w] \,, \hspace{5mm}
 [x_2,u_2] - [u_4,v]\,, \hspace{5mm}
[u_1,w] + a_2 [u_2,w] + b_2[u_3,w]  \,.
 \end{array}
 $$

Denote by $R\langle x_1, \ldots, x_n\rangle$ the
tensor algebra over a commutative ring $R$ with generators $x_1, \ldots ,
x_n$.   The image $A_X$ of $H_*(\Omega Z;\mathbb Z)$
in  $H_*(\Omega X;\mathbb Z)$ is  the quotient $H_*(\Omega Z;\mathbb
Z)/ I$, where $I$ denotes the two-sided ideal generated by the image of
$H_2(\Omega f)$ (\cite{Le}).
Thus $A_X = \mathbb Z\langle u_1, u_2, u_3, u_4, v, w, x_1, x_2\rangle /I$,
where all variables have degree one and $I$ is the ideal generated by the relations
above. For any field $\bk$, the Hilbert series of $A_X\otimes \bk$ and of $H_*(\Omega X;\bk)$
are related by a formula of Roos (\cite{JER}; Roos asserts this explicitly only for $\bk =\mathbb Q$,
but as Avramov observed in \cite{Av}
 it holds for any coefficient field)
$$P_{\Omega X}(t)^{-1}  = (1+t) A(t)^{-1} - t +8t^2 -13t^3\,.$$
In particular $H_*(\Omega X,\mathbb Z)$ and $A_X$ have the same torsion primes.

Our next step is to construct a graded algebra $E$ and then to identify the   algebra $A_X$   as
 the semi-tensor product
(\cite{Sm}) of the tensor algebra $\mathbb Z\langle x_1, x_2\rangle$
with   $E$.
Put $F = \mathbb Z \langle u_1, u_2, u_3, u_4, v,w\rangle$. For $i = 1, \ldots ,
4$,
define elements $\sigma_{i,m}, \rho_{i,m} \in F$ by $\sigma_{i,1} = u_i$, $\sigma_{i,m+1} =
[\sigma_{i,m},v]$ and   $\rho_{i,m} = [\sigma_{i,m-1},w]$.
Then
define    integers $a_{m}$  and $b_m$, $m\geq 3$,  by
$$\left(\begin{array}{l} a_{m}\\b_{m}\end{array}\right) = \left(\begin{array}{c} a\\
0\end{array}\right) + \left(\begin{array}{ll} b& c\\d&
0\end{array}\right)\left(\begin{array}{l}
a_{m-1}\\b_{m-1}\end{array}\right)\,.
$$
Finally, set
$$\tau_m = \rho_{1,m} + a_m \rho_{2,m} + b_m
\rho_{3,m}\,,\hspace{1cm} m\geq 2$$ and set
$$ S = \{\, \tau_m\,, m\geq 2\,\} \cup \{ \, a_m\rho_{4,m+1}\,,
m\geq 2\,\}\,.$$
Define $E$ to be the quotient of $F$ by the two-sided ideal $J$ generated
by $S$.

A unique action of $x_1$ and $x_2$ by derivations in $F$ is given
by
$$\left\{
\begin{array}{l}
x_1 *u_1 =[u_1,v] + a[u_2,v] \\
x_1*u_2 =  b[u_2,v] + d [u_3,v]\\
x_1 * u_3 = c[u_2,v]\\
x_1*u_4 = x_1*v=x_1*w = 0\\
x_2*u_1 = x_2*u_3 = x_2*u_4=x_2*v=x_2*w = 0\\
x_2*u_2 =  [u_4,v]\\
\end{array}
\right.$$
This extends uniquely to an action of $\mathbb Z\langle x_1, x_2
\rangle$ in $F$.
The action of $\mathbb Z\langle x_1, x_2\rangle$ preserves
$J$,  and so an action is induced in $F/J = E$.   Moreover, a presentation of the
 semi-tensor product of $\mathbb Z\langle x_1,
x_2\rangle $ with $E$ is given by
$$\mathbb Z\langle x_1, x_2 \rangle \otimes E  = \mathbb Z\langle u_1, u_2, u_3, u_4, x_1, x_2, v,w\rangle /
I\,,$$
where $I$ is the two-sided ideal generated by the relations
$[x_i,u_j]=x_i*u_j$, $[x_i,v]=x_i*v$, $[x_i,w]=x_i*w$ and
$  [u_1,w] + a_2 [u_2,w] + b_2[u_3,w] = 0
 $.
This identifies the semi-tensor product with $A_X$ and shows that
the torsion primes for $A_X$ are exactly those of $E$.

By construction, for $m\geq 3$, the elements $\rho_{4,m}$ are
torsion elements in $A_X$. In the same way as in \cite{An}, we have

\vspace{2mm}\noindent {\bf Lemma 1.} {\sl If $b$ is an integer, then $b\cdot \rho_{4,m}\neq 0$ in
$A_X$ unless $b\equiv 0$ mod $a_{m-1}$.}

\vspace{2mm}
Therefore all the prime numbers that appear as divisor of one of
the $a_{m}$ are torsion primes for $A_X$. The next lemma shows
that all the torsion primes are of that form.

\vspace{2mm}\noindent {\bf Lemma 2.}  {\sl If $p$ does not divide
one of the $a_m$, $m\geq 2$, then $A_X$ has no element of
order $p$.}

\vspace{2mm} \noindent {\bf proof:}  Invert the non-zero integers $a_m$, $m\geq 2$ and set
$R = \mathbb Z(\frac{1}{a_m})_{m\geq 2}$.
It is immediate from the criteria in \cite{Inert} that the elements $\rho_{4,m+1}$, $\tau_m$, $m\geq
2$,
form a strongly free (or inert) sequence in $F\otimes R$, and it follows that $E\otimes R = (F\otimes
R)/J$ splits back as a sub $R$-module of $F\otimes R$. In
particular, it is  torsion free. Hence so is $A_X\otimes R$.
 \hfill
$\Box$

\vspace{3mm}\noindent {\bf Proof of Theorem 1.} Consider the particular
case where $a=0, b = 4, c = -3, d= 1, a_2 = 11, b_2= 5$. In that
case
$$\tau_m = \rho_{1,m} + (2+3^m)\rho_{2,m} +
(2+3^{m-1})\rho_{3,m}$$
The integer $p$ is a homology torsion prime for $\Omega X$
if and only if for some $m\geq 2$, $2+3^m \equiv 0$ mod $p$.

We show now that the equation $2 + 3^m
\equiv
0$ mod $p$ has a solution when $p \equiv 5,7, 17$ or $19$ mod $24$, and does
not have a solution when $p \equiv 13$ or $23$ mod $24$.

In fact,  recall that a number a is a
 quadratic residue mod $p$ if and only if the Legendre symbol
$(a/p) = 1$. Recall (\cite{Bu}, p. 209, and Theorem 9-10),
that
$$(3/p) = \left\{ \begin{array}{lll} 1 & \mbox{if $p\equiv 1$ or
$11$}
&\mbox{mod 12}\\
-1 & \mbox{if $p\equiv$ 5 or 7}&\mbox{mod 12}\end{array}\right.
$$
$$(-2/p) = \left\{ \begin{array}{lll} 1 &\mbox{if $p\equiv 1$ or $3$}
&\mbox{mod $8$}\\
-1 & \mbox{if $p\equiv 5$ or $7$} &\mbox{mod $8$}
\end{array}
\right.
$$

Therefore if $p \equiv  5,7,17$ or $19$ mod 24, the integer $3$ is not a
 quadratic residue and is therefore a primitive root mod $p$.
Therefore the equation $2+3^m \equiv 0$ mod $p$ has always a solution.
On the other hand when $p\equiv 13$ or $23$ mod 24, then $3$ is a
quadratic residue, but $-2$ is not, and the equation $-2 \equiv 3^n$ mod $p$ has
therefore no solution.

Recall now   Dirichlet's Theorem (\cite{Bu}, Theorem 3-7): When $r$ and $s$ are relatively
prime integers  there are infinitely many primes
q such that $q \equiv r$ mod $s$. In particular there are infinitely many
prime numbers $p$ such that $p\equiv 5$ mod 24 and infinitely many
prime numbers such that
  $p \equiv 13$ mod 24. This shows that the set of torsion primes
for $A_X$ is an   infinite subset of $\cal P$ with an infinite complement.
\hfill $\square$

\vspace{3mm}\noindent {\bf Proof of Theorem 2.} Consider the particular
case where $b= 1, c = 0, d = 0, a_2=1, b_2=0$. In that case $a_m = 1 +
a(m-2)$. Thus $a_m\equiv 0$ mod $p$ for some $m$ if and only if $a\not\equiv 0$ mod $p$.
 Let   ${\cal S}= \{
p_1, \ldots , p_n\,\}$ be any finite subset of $\cal P$ and put
  $a = p_1 p_2 \cdots p_n$. Then $p$ is a prime divisor
for $A_X$ if and only if $p \not\in \cal S$. \hfill $\Box$

\vspace{5mm} Institut de Math\'ematiques,

Universit\'e Catholique de Louvain

1348 Louvain-La-Neuve, Belgium

\vspace{5mm} College of Computer, Mathematical and Physical Sciences

University of Maryland

College Park, MD 20742-3281, USA

\vspace{5mm} Facult\'e des Sciences

Universit\'e d'Angers

Bd. Lavoisier, 49045 Angers, France
\end{document}